\documentclass[pdflatex,sn-mathphys-num]{sn-jnl}
\usepackage[T1]{fontenc}
\usepackage{lmodern,amsmath,amssymb,amsthm,mathtools,microtype}
\newcommand{\E}{\mathbb E}
\newcommand{\Om}{\Omega}
\newcommand{\RR}{\mathbb R}
\newcommand{\Rtype}{T_p^{\mathrm R}}
\newcommand{\Etype}{T_p^{\mathrm E}}
\newcommand{\sign}{\operatorname{sign}}
\newcommand{\Lip}{\operatorname{Lip}}
\newcommand{\Id}{\operatorname{Id}}
\newcommand{\norm}[1]{\left\|#1\right\|}
\theoremstyle{thmstyleone}
\newtheorem{theorem}{Theorem}
\newtheorem{lemma}[theorem]{Lemma}
\newtheorem{proposition}[theorem]{Proposition}
\newtheorem{corollary}[theorem]{Corollary}
\theoremstyle{thmstylethree}
\newtheorem{remark}{Remark}
\newtheorem{example}{Example}
\hypersetup{pdftitle={Rademacher and Enflo types in Banach spaces: sharp comparison and orthogonality},
  pdfauthor={Zhiyao Fang and Qi Liu},
  pdfkeywords={Enflo type, Rademacher type, Birkhoff--James orthogonality, isosceles orthogonality, sharp constants}}

\begin{document}
\title[Rademacher and Enflo types in Banach spaces]{Rademacher and Enflo types in Banach spaces: sharp comparison and orthogonality}

\author[1]{\fnm{Zhiyao} \sur{Fang}}
\email{fzy73@stu.aqnu.edu.cn}

\author*[1]{\fnm{Qi} \sur{Liu}}
\email{liuq67@aqnu.edu.cn}
\affil[1]{%
  \orgdiv{School of Mathematics and Statistics},
  \orgname{Anqing Normal University},
  \orgaddress{\city{Anqing}, \postcode{246133}, \state{Anhui}, \country{China}}}

\abstract{Let \(1<p\le2\). We prove that the ratio of the Enflo to the Rademacher type \(p\) constant has supremum \(p\), already over finite-dimensional real or complex Banach spaces. This settles the sharpness question of Areshidze~\cite{aresh}, who improved the bound \(\pi/\sqrt{2}\) of Ivanisvili et al.~\cite{ivv} to \(p\). We also study orthogonal families in a real Banach space \(X\).
\begin{list}{}{\setlength{\leftmargin}{\dimexpr24pt+1.8em\relax}\setlength{\rightmargin}{24pt}\setlength{\labelwidth}{1.4em}\setlength{\labelsep}{0.4em}\setlength{\topsep}{2pt}\setlength{\itemsep}{0pt}}
\item[(a)] If \(p<2\), \(X\) has Enflo type \(p\) if and only if the Rademacher type \(p\) inequality holds uniformly over ordered families \((v_i)\) satisfying
\[
  v_i\perp_B v_j \qquad \text{for all } i<j,
\]
where \(\perp_B\) denotes Birkhoff--James orthogonality (that is, \(\|v_i+t v_j\|\ge\|v_i\|\) for every \(t\in\mathbb R\)).
\item[(b)] Uniform validity of the Rademacher type \(p\) inequality on pairwise isosceles orthogonal families implies that \(X\) has Enflo type \(s\) for every \(1\le s<p\), with a quadratic bound in the restricted type constant at \(s=2p/(p+1)\).
\end{list}
We also record the subcritical type implication for Pythagorean orthogonality and the sharp factor \(p\) for Roberts orthogonality. {For the exponents in (a), as the modulus of roundness \(\nu_X(p)\) approaches \(2\), the ratio of the Enflo to the Rademacher type \(p\) constant tends uniformly to \(1\) over real Banach spaces of Rademacher type \(p\). The same limit holds for the ratio of the two constants restricted by mutual Birkhoff--James orthogonality.}}
\maketitle

\section{Introduction}\label{sec:introduction}

Enflo type expresses a property of the metric of a Banach space through distances along the edges and antipodal pairs of a discrete cube \cite{enflo}. Rademacher type expresses the corresponding linear property through averages of signed sums. For background on type and cotype, see Hyt\"onen et al., \emph{Analysis in Banach Spaces}, Vol.~II \cite[pp.~53--162]{hvnwII}. The comparison belongs to the Ribe program \cite{naorRibe}, with antecedents in Enflo's work on uniform geometry \cite{enflo1969,enflo1970} and metric type \cite{bmw}. Mendel and Naor \cite[Theorem~1.1, p.~494]{mendelNaor} introduced scaled Enflo type and used discrete tori to give a metric characterization of Rademacher type \(p\) for \(1\le p\le2\). Their choice of scale preserves the type exponent. Ivanisvili et al.\ \cite[p.~666, Theorem~1.1]{ivv} later established the corresponding equivalence for ordinary Enflo type on the discrete cube, with a universal comparison factor of \(\pi/\sqrt{2}\). Pisier's cube inequality \cite{pisier1986} connects nonlinear type with vector-valued estimates on the cube \cite{naorSchechtman,hytonenNaor,eskenazisUMD,eskenazisNaor,talagrandCube,wagner}. For background on vector-valued martingales and Littlewood--Paley theory, see \cite{hvnwI}. {Related work on the discrete cube includes vector-valued influence and logarithmic Sobolev inequalities \cite{corderoEskenazisInfluence,corderoEskenazisLogSobolev,ivanisviliStone}. Concentration and heat-flow estimates on the biased cube also have applications to metric embeddings \cite{gordinBiasedCube,eskenazisHeatFlow}.}

Let \(\Om_n=\{-1,1\}^n\) with uniform probability. Products \(xy\) of elements of \(\Om_n\) are coordinatewise. Write \(x^j\) for \(x\) with coordinate \(j\) flipped and set \(D_jf(x)=(f(x)-f(x^j))/2\). Following Ivanisvili et al.\ \cite[Section~1]{ivv}, write \(\Rtype(X)\) and \(\Etype(X)\) for the least constants in the Rademacher and Enflo inequalities, respectively.
\begin{align}
 \big(\E\norm{\textstyle\sum_{j=1}^n s_jv_j}_X^p\big)^{1/p}
   &\le \Rtype(X)\big(\textstyle\sum_{j=1}^n\norm{v_j}_X^p\big)^{1/p},\label{eq:rademacher-type}\\
 \big\|(f(x)-f(-x))/2\big\|_{L_p(X)}
   &\le \Etype(X)\big(\textstyle\sum_{j=1}^n\norm{D_j f}_{L_p(X)}^p\big)^{1/p}.\label{eq:enflo-type}
\end{align}
The \(s_j\) are independent uniform signs. The quantifiers range over all \(n\), all \(v_j\in X\), and all \(f:\Om_n\to X\), respectively.

Ivanisvili et al.~\cite[Theorem~1.2 and Section~3]{ivv} obtain \(\pi/\sqrt{2}\) from the mass \(\pi/2\) of the heat-semigroup gradient kernel and a factor \(\sqrt{2}\) from the second-moment estimate after symmetrization. Areshidze's random-order martingale argument \cite[Theorem~1.1 and Corollary~3.1]{aresh} gives the upper bound \(p\) by summing weights proportional to \(m^{1/p-1}\), where \(m\) is the number of remaining coordinates. The same power law enters our lower bound through radial derivative profiles approaching \((1-s)^{1/p-1}\), whose integral on \([0,1]\) is \(p\). This settles Areshidze's sharpness question for \(1<p\le2\), already in finite dimensions.

\begin{theorem}\label{thm:sharp-comparison}
Fix \(1<p\le2\) and either scalar field \(\mathbb R\) or \(\mathbb C\). Then
\[
 C_p=\sup_X\frac{\Etype(X)}{\Rtype(X)}=p,
\]
where \(X\) ranges over nonzero Banach spaces of Rademacher type \(p\) over the chosen field. The same supremum is obtained over the finite-dimensional spaces over that field.
\end{theorem}

At \(p=1\), the triangle inequality for signed sums and the triangle inequality along cube paths give \(T_1^R(X)=T_1^E(X)=1\). For \(p>2\), repeating a unit vector \(n\) times would require \(\sqrt n\le T_p^R(X)n^{1/p}\), by comparison of the second and \(p\)-th moments. Thus \(1<p\le2\) covers the nontrivial range of the comparison.

Two notions of orthogonality extend the inner-product definition to real normed spaces. Birkhoff--James orthogonality, denoted by \(v\perp_B w\), means that \(\|v+tw\|\ge\|v\|\) for every \(t\in\mathbb R\), so \(v\) is a point of least norm on the line \(v+\mathbb Rw\) \cite{birkhoff1935,james}. Isosceles orthogonality, denoted by \(v\perp_I w\), means that \(\|v+w\|=\|v-w\|\) {\cite{james1945}}. Birkhoff--James orthogonality is homogeneous in each vector and can be asymmetric. Isosceles orthogonality is symmetric, and scaling one vector can change the relation. For a survey of both notions, see Alonso et al.\ \cite{alonso-martini-wu}. These notions have been studied in matrix distance problems, operator geometry, and norm-parallelism in Hilbert \(C^*\)-modules \cite{BhatiaSemrl1999,Bottazzi2019,MartinEtAl2024,ZamaniMoslehian2016}. Birkhoff--James orthogonality also characterizes the residuals of optimal polynomial approximation in spaces of analytic functions \cite{ChengFelder2023,ChengRossSeco2022}. Related directions include inner-product characterizations using Singer orthogonality and the study of isosceles orthogonal triples in linear \(2\)-normed spaces \cite{Alonso1993,Cho1992}.

Liu et al.\ obtained exact representations of four unit-sphere geometric constants, including the modified von Neumann--Jordan constant, using isosceles orthogonal pairs and, through a norming-functional construction, Birkhoff--James orthogonal pairs \cite{liu-wang-bao-geometric}. This raises a related question for inequalities involving arbitrarily long vector families. Here we ask when signed-sum estimates on orthogonal families control Enflo type on the whole space and whether orthogonality restrictions preserve the optimal comparison factor.

For a real Banach space \(X\), define the restricted type constants as follows.
\begin{itemize}
\item \(B_p(X)\) is the best constant in \eqref{eq:rademacher-type} over families \((v_i)\) in \(X\) with \(v_i\perp_B v_j\) for every ordered pair \(i\ne j\).
\item \(I_p(X)\) is the best constant in \eqref{eq:rademacher-type} over families \((v_i)\) in \(X\) with \(v_i\perp_I v_j\), that is, \(\|v_i+v_j\|=\|v_i-v_j\|\), for every pair \(i<j\).
\end{itemize}

For comparison, the unrestricted type constants satisfy
\begin{align*}
 \Rtype(X)\le \Etype(X)&\le \frac{\pi}{\sqrt{2}}\Rtype(X)
 &&\text{\cite[p.~666, Theorem~1.1]{ivv}},\\
 \Etype(X)&\le p\Rtype(X)
 &&\text{\cite[p.~5, Corollary~3.1]{aresh}}.
\end{align*}
Here the Rademacher inequality is taken over all finite vector families. The next two theorems show what follows when the Rademacher inequality is imposed only on the orthogonal families defined above.

\begin{theorem}\label{thm:orthogonal-tests}
For every \(1<p<2\), there is a finite constant \(K_p\) depending only on \(p\) such that every nonzero real Banach space \(X\) satisfies
\[
 B_p(X)\le T_p^R(X)\le K_p B_p(X)^2,
 \qquad
 T_p^E(X)\le pK_p B_p(X)^2.
\]
In particular, \(X\) has ordinary Enflo type \(p\) if and only if \(B_p(X)<\infty\).
\end{theorem}

The same equivalence holds for ordered families with \(v_i\perp_B v_j\) for \(i<j\), since their best constant lies between \(B_p(X)\) and \(T_p^R(X)\).\footnote{At \(p=2\), the entropy estimate in Lemma~\ref{lem:unit-gaussian} incurs a logarithmic loss, and the dyadic error bound in Lemma~\ref{lem:gaussian-transfer} loses its geometric decay. These estimates therefore do not give a dimension-independent bound for the reverse implication at the endpoint.}

\begin{theorem}\label{thm:isosceles-tests}
For \(1<p\le2\), put \(r=2p/(p+1)\). There is a finite constant \(K_p^{I}\) depending only on \(p\) such that every nonzero real Banach space \(X\) satisfies
\[
 T_r^R(X)\le K_p^{I} I_p(X)^2,\qquad
 T_r^E(X)\le rK_p^{I} I_p(X)^2.
\]
In particular, if \(I_p(X)<\infty\), then \(X\) has Enflo type \(r\).
\end{theorem}

Theorems~\ref{thm:orthogonal-tests} and~\ref{thm:isosceles-tests} connect orthogonality to the Ribe program by deriving Enflo type from Rademacher inequalities on Birkhoff--James and isosceles orthogonal families. Theorem~\ref{thm:nj-stability} links the stability of the comparison between linear and metric type constants to a four-point distance inequality.

The linear cube map \(f(x)=\sum_j x_jv_j\) gives \(\Rtype(X)\le\Etype(X)\). The reverse bound \(\Etype(X)\le p\Rtype(X)\) for \(1\le p\le2\) is proved in \cite[p.~5, Corollary~3.1]{aresh}. Areshidze's proof of Theorem~1.1 \cite[p.~4]{aresh} starts from the martingale telescoping identity
\[
 f-\E f=M_n^\pi-M_0^\pi=\sum_{k=1}^n d_{\pi_k}^{A_{k-1}^\pi} f,
\]
where \(\pi\) is a permutation of the coordinates, \(A_{k-1}^\pi=\{\pi_1,\ldots,\pi_{k-1}\}\), \(M_k^\pi\) is the conditional mean of \(f\) given the first \(k\) coordinates in this order, and \(d_{\pi_k}^{A_{k-1}^\pi}f=M_k^\pi-M_{k-1}^\pi\). Areshidze averages this identity over \(\pi\). For the orthogonality-restricted inequalities in Section~\ref{sec:orthogonal-tests}, we use the random-sign representation below, in which Bernoulli selectors retain a subfamily of the edge vectors. Put \(P_{\rm odd}f(x)=(f(x)-f(-x))/2\) and \(H_p(f)=(\sum_j\|D_jf\|_{L_p(X)}^p)^{1/p}\). For \(0<u<1\), let \(I_j\) be independent Bernoulli variables of mean \(1-u\), independent of the signs \(s_j\), and put \(\eta_j=s_j\) if \(I_j=1\) and \(\eta_j=1\) otherwise. Then
\begin{equation}\label{eq:random-sign-representation}
 f(x)-\E f=\int_0^1 \E_{I,s}\sum_j\frac{I_j s_j}{1-u}\,D_j f(x\eta)\,du .
\end{equation}
For the Walsh terminology, see O'Donnell's \emph{Analysis of Boolean Functions} \cite[\S\S1.2--1.4 and \S2.2, 2021 edition]{odonnell-abf}. For \(S\subseteq[n]\), the Walsh character on \(\Om_n\) is \(\chi_S(x)=\prod_{i\in S}x_i\), with \(\chi_\varnothing=1\) \cite[pp.~22--23, \S1.2 and Definition~1.2, 2021 edition]{odonnell-abf}. Its Walsh degree is \(|S|\) \cite[p.~28, text before Definition~1.19]{odonnell-abf}. The Walsh coefficients are \(\widehat f(S)=\E_x[f(x)\chi_S(x)]\), so \(f=\sum_{S\subseteq[n]}\widehat f(S)\chi_S\) \cite[p.~25, Proposition~1.8, 2021 edition]{odonnell-abf}. The Walsh degree of a nonzero \(f\) is \(\max\{|S|:\widehat f(S)\ne0\}\). On a Walsh character of degree \(d>0\), the integrand has multiplier \(du^{d-1}\), whose integral is one. Both sides vanish on constants. By Jensen's inequality and the substitution \(z=x\eta\), the Rademacher inequality bounds the integrand in \(L_p(X)\) by \((1-u)^{1/p-1}\Rtype(X)H_p(f)\), using \(\E I_j=1-u\). Integration and the contractivity of \(P_{\rm odd}\) give \(\Etype(X)\le p\Rtype(X)\).

Sections~\ref{sec:radial} and~\ref{sec:sharpness} prove sharpness using radial gradients and an operator decomposition. Section~\ref{sec:orthogonal-tests} proves the main orthogonality results via Gaussian blocks and records the corresponding conclusions for Pythagorean and Roberts orthogonality. {Section~\ref{sec:nj-stability} proves stability of the type comparison as the modulus of roundness tends to two.}

\section{Rademacher sums of radial gradients}\label{sec:radial}

For the lower bound we first work over \(\RR\) and set \(q=p/(p-1)\). Give \([n]\) counting measure. Let \(R_m:\ell_p^m\to L_p(\Om_m)\) be the map \(a\mapsto\sum_j s_j a_j\). For \(U:V\to Z\), set \(R_p(U)=\sup_m\|R_m\otimes U\|\), where the \(m\)-th operator sends \((v_j)_{j=1}^m\in\ell_p^m(V)\) to the function \(s\mapsto\sum_j s_jUv_j\) in \(L_p(\Om_m;Z)\). Thus \(R_p(\Id_X)=\Rtype(X)\). The scalar bound \(\|R_ma\|_p\le\|R_ma\|_2=\|a\|_2\le\|a\|_p\) is attained on coordinate vectors.

For an even nonnegative Lipschitz function \(\psi:[-1,1]\to\mathbb R\) satisfying
\begin{equation}\label{eq:profile-bound}
 \psi(t)\le(1-|t|)^{-1/q}\qquad (|t|<1),
\end{equation}
set \(\Phi(t)=\int_0^t\psi(v)\,dv\) and \(M(x)=n^{-1}\sum_{i=1}^n x_i\). For \(z\in\Om_n\) define \(f_z(x)=n\Phi(M(xz))\). The linear map
\[
 U_n:\ell_1(\Om_n\times[n])\longrightarrow\ell_\infty(\Om_n),
 \qquad (U_ne_{x,j})(z)=D_jf_z(x),
\]
records the gradients of these translates. Here \(e_{x,j}\) denotes the corresponding coordinate vector. We will show that \(R_p(U_n)\le n^{1-1/p}(1+o(1))\). The following estimate is uniform in the number of terms and allows repeated vertices and coordinates.

\begin{lemma}\label{mask}
Fix a nonnegative Lipschitz function \(\psi:[0,1]\to\RR\) satisfying \(\psi(t)\le(1-t)^{-1/q}\) for \(t<1\). Given \(n\), choose any \(m\), numbers \(b_\ell\ge0\) with \(\sum_{\ell=1}^m b_\ell^p=n\), indices \(j_\ell\in[n]\), and \(h_\ell\in\Om_n\). For independent signs \(s_\ell\) put \(m_\ell(\eta)=n^{-1}\langle h_\ell,\eta\rangle\). Uniformly in all these choices,
\[
 \left\|\sup_{\eta\in\Om_n}\left|
 \frac1n\sum_{\ell=1}^m b_\ell s_\ell\eta_{j_\ell}h_{\ell,j_\ell}
     \psi(|m_\ell(\eta)|)\right|\right\|_{L_p(s)}
 \le 1+o_{p,\psi}(1)\qquad(n\to\infty).
\]
\end{lemma}
\begin{proof}
Select \(\eta=\eta(s)\) attaining the absolute maximum, breaking ties by a fixed ordering of \(\Om_n\). The argument applies to any such sign-dependent choice. We may suppose \(K=\|\psi\|_\infty>0\).

A centered variable \(Z\) with \(\E e^{\lambda Z}\le e^{\lambda^2 v^2/2}\) for all real \(\lambda\) is subgaussian with variance proxy \(v^2\) \cite[\S2.3]{boucheron-lugosi-massart}. For \(N\) such variables with a common proxy \(v^2\), the \(L_p\) norm of the maximum of their absolute values is at most \(C_p'v\sqrt{\log(2N)}\). This follows by integrating the union bound \(2N e^{-u^2/(2v^2)}\) above \(u=v\sqrt{2\log(2N)}\). A signed linear sum has this bound with \(v\) equal to the Euclidean norm of its coefficients, by independence and \(\cosh u\le e^{u^2/2}\).

\emph{Approximating the row averages.} Fix \(M\ge1\). Discard the terms \(b_\ell>M\) at a pointwise cost \(\le K M^{1-p}\), because their total weight is at most \(M^{1-p}n\). In what follows all row sums are over the retained terms. When \(B^2=\sum b_\ell^2\) is positive, we have \(B^2\le M^{2-p}n\). If \(B=0\), all retained terms vanish. The matrix \(H\) with rows \(b_\ell h_\ell/(B\sqrt n)\) has Hilbert--Schmidt norm one, since the sum of the squares of its entries is one. Put \(r=\lfloor\sqrt n\rfloor\) and use the eigenvectors of \(H^*H\) for its \(r\) largest eigenvalues, including zero eigenvalues if necessary. Partition \(\Om_n\) by quantizing the coordinates of \(\eta/\sqrt n\) in those directions to intervals of length \(1/r\). The partition has \(N\le(2r+2)^r\) cells and \(\log(2N)=o(n)\).

Let \(T\) be the cell of the selected \(\eta\). For each cell \(t\) of positive probability set \(d^t=\E(\eta\mid T=t)\), \(\mu_\ell^t=n^{-1}\langle h_\ell,d^t\rangle\), and \(\theta_t=n^{-1}\sum_j|d_j^t|\). For every \(\eta\) in this cell,
\begin{equation}\label{svd}
 \sum_\ell b_\ell^2|m_\ell(\eta)-\mu_\ell^t|^2
 \le B^2 e_n^2,\qquad
 e_n=r^{-1/2}+2(r+1)^{-1/2}.
\end{equation}
Indeed, the left side divided by \(B^2\) is \(\|H(\eta-d^t)/\sqrt n\|_2^2\). The first \(r\) coordinates of the input have moduli at most \(1/r\). The remaining input has Euclidean norm at most two, and the \((r+1)\)-st singular value is at most \(1/\sqrt{r+1}\) (zero if \(r=n\)). The operator norm of \(H\) is at most one. This proves \eqref{svd}. For each fixed \(\eta\), the signed error obtained by replacing \(\psi(|m_\ell(\eta)|)\) by \(\psi(|\mu_\ell^t|)\) has variance proxy at most \((B e_n\Lip(\psi)/n)^2\). Each \(\eta\) is tested with its own cell's deterministic vector \(d^t\). Maximizing over at most \(2^n\) deterministic choices gives an \(L_p\) error \(O_p(\Lip(\psi)e_n B/\sqrt n)=o(1)\) for fixed \(M\).

\emph{Concentration within the coordinate groups.} Since \(|\mu_\ell^t|\le\theta_t\), set
\[
 W_t=\min\{K,(1-\theta_t)^{-1/q}\},\qquad
 c_\ell^t=h_{\ell,j_\ell}\psi(|\mu_\ell^t|)/W_t ,
\]
where \(W_t=K\) for \(\theta_t=1\). Thus \(|c_\ell^t|\le1\). Put
\[
 Z_j^t=\sum_{\ell:j_\ell=j}b_\ell c_\ell^t s_\ell,\qquad
 V_j^2=\sum_{\ell:j_\ell=j}b_\ell^2 .
\]
For each fixed \(t\), the \(Z_j^t\) are independent and symmetric under the original sign law, and each has subgaussian variance proxy \(V_j^2\). Also \(\sum V_j^2=B^2\) and \(\sum V_j^p\le\sum b_\ell^p\le n\). We estimate simultaneous maxima over cells under the original sign law, allowing \(T\) to depend on all signs.

Fix \(L\ge1\) and discard groups \(V_j>L\). For each \(t\) the mean of \(\sum_{j:V_j>L}|Z_j^t|\) is at most \(\sum_{V_j>L} V_j\le nL^{1-p}\). Changing \(s_\ell\) changes this sum by at most \(2b_\ell\). Successively exposing the uniform signs yields conditional symmetric martingale increments of absolute value at most \(b_\ell\). Iterating their exponential moment bounds makes the centered sum subgaussian with variance proxy \(B^2\). The \(L_p\) error in \(W_Tn^{-1}\sum_j\eta_j Z_j^T\) from discarding the groups is consequently at most
\[
 K\left(L^{1-p}+C_p'B\sqrt{\log(2N)}/n\right).
\]
Let \(J=\{j:V_j\le L\}\), \(z_j^t=\sign(d_j^t)\), choosing \(+1\) at zero. A second application of the maximal bound gives
\[
 \left\|\max_t \left|n^{-1}\sum_{j\in J} z_j^t Z_j^t\right|\right\|_{L_p}
       \le C_p' B\sqrt{\log(2N)}/n=o(1).
\]
Moreover, with \(u_+=\max(u,0)\),
\begin{equation}\label{exc}
 H_n=\left(\max_{t,\sigma=\pm1}\frac1n
    \sum_{j\in J}(\sigma z_j^t Z_j^t)_+^p-\frac12\right)_+
       \quad\text{satisfies}\quad \E H_n=o(1).
\end{equation}
In fact, the independent summands for a fixed \(t,\sigma\) have means bounded by \(V_j^p/2\), by symmetry and the \(L_2\) bound. Their total mean is \(\le n/2\). Since \(\Pr(|Z_j^t|>v)\le 2e^{-v^2/(2L^2)}\) and \(p\le2\), their centered versions \(Y_j\) have uniformly bounded absolute exponential moments at some positive parameter depending only on \(p,L\). The inequality \(|e^v-1-v|\le v^2e^{|v|}\) gives \(\E e^{\lambda Y_j}\le e^{C_{p,L}\lambda^2}\) for \(|\lambda|\le c_{p,L}\). Independence bounds the moment generating function of each centered sum by \(e^{n C_{p,L}\lambda^2}\). Bounding the exponential of the maximum by the sum of exponentials (including the zero variable) and taking \(\lambda=\sqrt{\log(2N+1)/n}\) for large \(n\) proves \eqref{exc}.

\emph{Using the conditional mean of the selected signs.} Put \(u=n^{-1}\#\{j:\eta_j\ne z_j^T\}\). By the definition of \(d^t\),
\(\E(u\mid T=t)=(1-\theta_t)/2\). The identity \(\eta_j=z_j^T(1-2\mathbf1_{\{\eta_j\ne z_j^T\}})\) separates the centered signed sum just bounded. For the remaining sum, H\"older's inequality for each sign \(\sigma\) yields
\[
 \sigma\, n^{-1}\sum_{j\in J}z_j^T Z_j^T
 \mathbf1_{\{\eta_j\ne z_j^T\}}
 \le (1/2+H_n)^{1/p}u^{1/q}
 \le 2^{-1/p}u^{1/q}+H_n^{1/p}.
\]
The term \(2W_TH_n^{1/p}\) tends to zero in \(L_p\) by \eqref{exc}. As \(p/q=p-1\le1\), Jensen conditional on \(T\) gives
\[
 \norm{2^{1-1/p}W_T u^{1/q}}_{L_p}^p
 \le \E W_T^p(1-\theta_T)^{p/q}\le1 .
\]
This also holds on \(\theta_T=1\), where \(u=0\) almost surely. Combining the estimates at fixed \(M,L\) bounds the limit superior by \(1+K(M^{1-p}+L^{1-p})\), uniformly in all arrays. Let \(M,L\) tend to infinity.
\end{proof}

\begin{proposition}\label{radtype}
For fixed \(p\) and an even nonnegative Lipschitz \(\psi\) satisfying \eqref{eq:profile-bound}, the associated operator satisfies \(R_p(U_n)\le n^{1-1/p}(1+o(1))\) as \(n\to\infty\).
\end{proposition}
\begin{proof}
For any number of inputs and any fixed lengths \(a_\ell\ge0\), the function \(\E\|\sum_\ell s_\ell a_\ell U_n v_\ell\|_\infty^p\) is separately convex in \(v_\ell\). The unit ball of \(\ell_1(\Om_n\times[n])\) is the convex hull of signed coordinate vectors, so it suffices to test tuples of those vectors. Absorb the signs into \(s_\ell\) and normalize \(\sum_\ell a_\ell^p=n\). For the selected vertices \(h_\ell\) and indices \(j_\ell\), writing \(m_\ell(z)=\langle z,h_\ell\rangle/n\), we have
\[
 (U_n e_{h_\ell,j_\ell})(z)=h_{\ell,j_\ell}z_{j_\ell}\psi(m_\ell(z))+r_\ell(z),
 \qquad |r_\ell(z)|\le 2\Lip(\psi)/n .
\]
Indeed, dividing the difference of \(\Phi\) by the increment \(2h_{\ell,j_\ell}z_{j_\ell}/n\) gives the average of its derivative over that interval. By the maximal bound over \(2^n\) choices of \(z\), the \(L_p\) norm of \(\sup_z n^{-1}|\sum s_\ell a_\ell r_\ell(z)|\) is at most \(C_p'\Lip(\psi)n^{1/p-3/2}\), using \(\|a\|_2\le\|a\|_p=n^{1/p}\). Lemma~\ref{mask} bounds the principal sum divided by \(n\) by \(1+o(1)\). Dividing the resulting bound by \(n^{1/p}\) proves the claim.
\end{proof}

\section{The sharp comparison constant}\label{sec:sharpness}

We now prove the lower bound in Theorem~\ref{thm:sharp-comparison}. If a factor \(\kappa<p\) sufficed, the operator bipolar theorem would give a decomposition of the inverse cube gradient on each fixed range of odd Walsh degrees. We test this decomposition on the radial gradients from Section~\ref{sec:radial} and show that it forces \(\kappa\ge p\). The resulting finite-dimensional real examples then give the complex lower bound by an isometric complexification that preserves the Rademacher type constant.

Write \(Df=(D_jf)_{j\in[n]}\) and define \(\Gamma_n(Df)=f\) on gradients of odd scalar functions in \(L_p(\Om_n\times[n])\). The Walsh gradient \(D\chi_S\) is the tuple of coordinate derivatives of \(\chi_S\). Its \(j\)-th component is \(D_j\chi_S=\mathbf{1}_{\{j\in S\}}\chi_S\), since flipping \(x_j\) changes the sign precisely when \(j\in S\) (compare the discrete derivative convention in O'Donnell \cite[p.~47, Definition~2.16 and Proposition~2.19, 2021 edition]{odonnell-abf}). A function with zero gradient is constant, so \(\Gamma_n\) is well-defined. Odd Walsh characters span the odd functions, also with coefficients in \(X\). Hence \(\|\Gamma_n\otimes\Id_X\|\le\Etype(X)\) by \eqref{eq:enflo-type}.

The operator bipolar theorem gives an approximate decomposition of the inverse gradient from a bound on its Banach-valued amplifications. We use the decomposition on each fixed range of Walsh degrees.

For an operator \(S\) on a subspace of scalar \(L_p\), with values in another \(L_p\), the regular norm is \(\sup_k\|S\otimes \Id_{\ell_\infty^k}\|\), using the ambient vector-valued norms. For the regular norm on operators defined on subspaces of Banach lattices, see Pisier \cite[formula~(1) and the discussion preceding Theorem~3 in the preprint version]{pis}. A finite full matrix \(S\) has regular norm \(\||S|\|_{p\to p}\), where \(|S|\) replaces each coefficient by its absolute value. The upper bound follows from the triangle inequality on each row. For the reverse bound test on a positive input, adjoining one \(\ell_\infty\) coordinate for each row and signing the input in that coordinate according to that row. Thus every Banach amplification of a full finite matrix is bounded by its regular norm.

\begin{proposition}\label{prop:decomposition}
Fix \(\kappa\ge1\). Let \(\mathcal N\subset\mathbb N\) be infinite and \(c_n>1\) with \(c_n\to\infty\) along \(\mathcal N\). Suppose that for every \(n\in\mathcal N\) and every real Banach space \(X\),
\begin{equation}\label{eq:localized-implication}
 T_p^R(X)\le c_n \quad\Longrightarrow\quad
 \|\Gamma_n\otimes\Id_X\|\le\kappa c_n .
\end{equation}
For each \(n\in\mathcal N\) there exist full finite regular contractions \(A,B,D_0\) and a finite \(\ell_p\)-direct sum \(R:P\to Q\) of the \(R_m\) such that
\[
 A:L_p(\Om_n\times[n])\to L_p(\Om_n),\qquad
 D_0:L_p(\Om_n\times[n])\to P,\qquad B:Q\to L_p(\Om_n).
\]
The operator \(E_n=\Gamma_n-\kappa c_n A-\kappa BRD_0\) on odd gradients satisfies
\begin{equation}\label{scalarerror}
 \max_{\substack{S\subset[n],\ |S|\le h\\ |S|\ \mathrm{odd}}}
 \|E_n(D\chi_S)\|_p=o(1)\qquad(n\to\infty,\ n\in\mathcal N)
\end{equation}
for each fixed positive integer \(h\).
\end{proposition}

\begin{proof}
Fix \(n\in\mathcal N\) and put \(c=c_n\). If a real Banach space \(X\) satisfies \(\sup_m\|(R_m/c)\otimes\Id_X\|\le1\), then \(T_p^R(X)\le c\), so \eqref{eq:localized-implication} gives \(\|(\Gamma_n/(\kappa c))\otimes\Id_X\|\le1\).

\emph{Application of the bipolar theorem.}
Apply de la Salle's operator bipolar theorem \cite[Definition~1.5 and Theorem~1.6]{dls} in the real setting to this implication. Use the entire finite list of Walsh gradients \(D\chi_S\) with \(|S|\) odd, which spans the domain of \(\Gamma_n\). Write
\[
 G=L_p(\Om_n\times[n]\times[0,1]),\quad
 Y=L_p(\Om_n\times[0,1]),\quad H=L_p([0,1]).
\]
Use Lebesgue probability on the new factors and identify original arrays with their extensions constant in that variable. For every \(\epsilon>0\) the theorem supplies a finite \(\ell_p\)-direct sum \(R:P\to Q\) of operators \(R_m\), an operator \(S_0\) with regular norm at most one on a subspace of \(L_p\), and spatial isometries \(U,V\) on subspaces such that, for each test \(S\), there are vectors with
\begin{align}
 U(g_S,h_S)=(u_S,a_S),\qquad
 &V(S_0u_S,Ra_S/c)=(g'_S,h_S),\label{bip}\\
 \norm{g_S-D\chi_S}_{G}\le\epsilon,\qquad
 &\norm{g'_S-\chi_S/(\kappa c)}_Y\le\epsilon .\notag
\end{align}
Here \(g_S\in G,\ h_S\in H,\ g'_S\in Y,\ a_S\in P\), and \(u_S\) lies in the domain of \(S_0\). A spatial isometry is a composition of changes of density and maps that preserve the joint distributions of finite tuples away from the origin. Each such isometry preserves the norms of finite \(\ell_\infty\) amplifications. A change of density preserves the integral of the \(p\)-th power of the maximum modulus by cancellation of the density factors. An equidistribution map preserves this integral by equality of the joint distributions away from the origin. Thus \(U,V\) are regular contractions.

Pisier's regular extension theorem \cite[Theorem~3 in the preprint version]{pis} extends each of \(U,V,S_0\) to its full ambient domain with values in its original ambient target and with the same regular norm. The domain can first be closed by continuity. To obtain the real version of the extension, complexify a real subspace map \(C\) to \(C_{\mathbb C}\). For complex inputs \(v_j\), the pointwise maximum \(\max_j|C_{\mathbb C}v_j|\) equals the supremum of \(|C\operatorname{Re}(e^{i\theta}v_j)|\) over \(j\) and a countable dense set of angles. Applying real regularity to finite sets of angles and using monotone convergence shows that complexification preserves the regular norm. Restrict the complex extension to real inputs and take the real part to obtain the real extension.

\emph{Eliminating the auxiliary coordinate.}
Using these full extensions, form a linear map
\[
 \mathcal L:(G\oplus_p Q)\oplus_p H\longrightarrow (Y\oplus_p P)\oplus_p H .
\]
Given \((g,b,h)\) compute \((u,a)=U(g,h)\), then \((y,h')=V(S_0u,b)\), and return \((y,a,h')\). In any finite \(\ell_\infty\) amplification,
\[
 \|y\|^p+\|h'\|^p+\|a\|^p
 \le\|S_0u\|^p+\|b\|^p+\|a\|^p
 \le\|g\|^p+\|h\|^p+\|b\|^p.
\]
Consequently \(\mathcal L\) is a regular contraction. By \eqref{bip}, if \(b_S=Ra_S/c\), then \(\mathcal L(g_S,b_S,h_S)=(g'_S,a_S,h_S)\).

Write \(\mathcal L=\left(\begin{smallmatrix}J&K\\ M_0&N\end{smallmatrix}\right)\) with the auxiliary \(H\) last. For an external input \(x\in G\oplus_p Q\), define \(h_r=r(I-rN)^{-1}M_0x\) and \(V_r x=Jx+Kh_r\) for \(0<r<1\). The Neumann series \((I-rN)^{-1}=\sum_{a\ge0}r^aN^a\) converges in each finite \(\ell_\infty\) amplification since \(N\) is a regular contraction. The second output at \((x,h_r)\) equals \(h_r/r\), so
\[
 \|V_r x\|^p+\|h_r/r\|^p\le \|x\|^p+\|h_r\|^p
\]
in each amplification. Thus \(V_r\) is a regular contraction. If \(h=M_0x+Nh\), then \(v_r=h-h_r=(1-r)(I-rN)^{-1}h\) satisfies \(\|v_r\|\le\|h\|\) and
\[
 \|Nv_r-v_r\|\le(1-r)(\|v_r\|+\|h\|)/r\longrightarrow0 .
\]
Contractivity on \((0,v_r)\) gives \(\|Kv_r\|^p\le\|v_r\|^p-\|Nv_r\|^p\to0\), using boundedness. Hence \(V_rx\to Jx+Kh\) whenever such \(h\) exists. Given \(\delta>0\), finiteness of the test list lets us choose a common \(r\) for which every \((g_S,b_S)\) has output within \(\delta\) of \((g'_S,a_S)\).

In external blocks write \(V_r=\left(\begin{smallmatrix}A&B\\ D_0&F_0\end{smallmatrix}\right)\). Its four blocks are full regular contractions. We have
\[
 \|g'_S-Ag_S-Bb_S\|\le\delta,\qquad
 \|a_S-D_0g_S-F_0b_S\|\le\delta .
\]
Scalar \(\|R\|\le1\). With \(b_S=Ra_S/c\), the second inequality gives
\(\|a_S\|\le(\|g_S\|+\delta)/(1-1/c)\). Substitution into the first yields
\begin{equation}\label{smallremainder}
 \|g'_S-Ag_S-c^{-1}BRD_0g_S\|
 \le \delta(1+1/c)+\frac{\|g_S\|+\delta}{c(c-1)} .
\end{equation}

Insert constant functions along \([0,1]\) at the edge input and integrate out \([0,1]\) at the cube output. These maps are regular contractions, so the compressed \(A,B,D_0\), still with these names, are full finite matrices of regular norm at most one. Put \(\epsilon=\delta=c^{-2}n^{-2}\). Using \(\|D\chi_S\|_p=|S|^{1/p}\), multiply \eqref{smallremainder} by \(\kappa c\) and use \eqref{bip}. We have obtained a single triple \(A,B,D_0\) for the cube such that the operator on scalar odd gradients
\[
 E_n=\Gamma_n-\kappa c A-\kappa BRD_0
\]
satisfies
\[
 \sup_{\substack{S\subset[n],\ |S|\le h\\ |S|\ {\rm odd}}}
 \|E_n(D\chi_S)\|_p
 \le \frac{\kappa h^{1/p}}{c-1}+o(1).
\]
Indeed, replacing the two perturbed tests in \eqref{smallremainder} costs at most \(\kappa c\epsilon(2+1/c)\), besides the terms with \(\delta\) and the additive \(\epsilon\) in \(\|g_S\|\). Since \(c_n\to\infty\), this proves \eqref{scalarerror}.
\end{proof}

\begin{lemma}\label{lem:localized-witness}
Fix \(1\le\kappa<p\) and put \(c_n=2+n^{(1-1/p)/2}\). For every sufficiently large \(n\), there are a finite-dimensional real space \(X_n\) and a map \(f:\Om_n\to X_n\) with \(H_p(f)>0\) such that
\begin{equation}\label{eq:localized-witness}
 T_p^R(X_n)\le c_n,\qquad
 \|P_{\rm odd}f\|_{L_p(X_n)}>\kappa c_n H_p(f).
\end{equation}
\end{lemma}

\begin{proof}
If the assertion fails along an infinite set \(\mathcal N\), then \eqref{eq:localized-implication} holds along that set. Indeed, a violation for an arbitrary real space would be witnessed by one odd cube map, whose values span a finite-dimensional space with Rademacher type constant at most that of the original space. Apply Proposition~\ref{prop:decomposition} along \(\mathcal N\), where \(c_n\to\infty\) and \(c_n n^{1/p-1}\to0\).

Fix an even nonnegative polynomial \(\psi\) satisfying \eqref{eq:profile-bound}. Use \(f_z(x)=n\Phi(M(xz))\) from Section~\ref{sec:radial}, and put \(F(x)=\Phi(M(x))\). If \(\Phi(t)=\sum_{k=0}^h a_k t^k\), the Walsh coefficients of \(F\) satisfy
\[
 \widehat F(S)=0 \text{ unless } |S|\le h\text{ is odd},
 \qquad \sum_S |\widehat F(S)|\le\sum_{k=0}^h|a_k| .
\]
Indeed \(M\) has nonnegative Walsh coefficients summing to one, as does \(M^k\) by multiplication of characters. Oddness removes the even degrees. Thus testing \eqref{scalarerror} via the expansion \(f_z(x)=n\sum_S \widehat F(S)\chi_S(z)\chi_S(x)\) gives
\[
 \left|n^{-1}\E_z [E_n(Df_z)](z)\right|
 \le\sum_S |\widehat F(S)|\,\|E_n(D\chi_S)\|_p=o(1).
\]
The contribution of \(\Gamma_n\) to this normalized diagonal average equals \(\Phi(1)\). For \(A\), \(|D_j f_z(x)|\le\|\psi\|_\infty\) by the mean-value theorem. The input constant array has norm \(n^{1/p}\) and \(|A|\) is contractive, so
\[
 n^{-1}\left|\E_z [A(Df_z)](z)\right|
 \le n^{-1}\|\psi\|_\infty \E_z(|A|\mathbf1)(z)
 \le\|\psi\|_\infty n^{1/p-1}.
\]
Hence even after multiplication by \(c\) its contribution tends to zero.

Finally \(D_jf_z(x)=(U_ne_{x,j})(z)\). Write \(V_n=\ell_1(\Om_n\times[n])\) and let \(g(x,j)=e_{x,j}\). Its \(L_p(\Om_n\times[n];V_n)\) norm is \(n^{1/p}\). At output row \(z\), evaluating the \(\ell_\infty(\Om_n)\)-valued array at coordinate \(z\) gives the diagonal entry. In particular,
\begin{align*}
 n^{-1}\big|\E_z[BRD_0(Df_z)](z)\big|
 &\le n^{-1}\|(B\otimes\Id_{\ell_\infty})(R\otimes U_n)
             (D_0\otimes\Id_{V_n})g\|_{L_p(\ell_\infty)}\\
 &\le n^{1/p-1}R_p(U_n)\le 1+o(1).
\end{align*}
Here \(R\otimes U_n\) has norm at most \(R_p(U_n)\), by summing the \(p\)-th powers over the parallel blocks. Proposition~\ref{radtype} gives the last inequality. The three contributions to the diagonal test yield
\(\Phi(1)=\int_0^1\psi(t)\,dt\le\kappa\).

For \(L\ge1\), extend \(v_L(t)=\min\{L,(1-|t|)^{-1/q}\}\) continuously to \([-1,1]\). Given \(0<\alpha<1/2\), approximate \(v_L\) uniformly to within \(\alpha\) by an even real polynomial \(P\), using Weierstrass approximation and symmetrization. Then \(0\le P-\alpha\le v_L\), since \(v_L\ge1\), and its integral on \([0,1]\) differs from that of \(v_L\) by at most \(2\alpha\). Taking \(\psi=P-\alpha\), then letting \(\alpha\downarrow0\) and \(L\uparrow\infty\), gives
\[
 \kappa\ge\int_0^1(1-t)^{-1/q}\,dt=p.
\]
This contradicts \(\kappa<p\).
\end{proof}

To transfer the lower bound from real to complex spaces, we choose a complexification that preserves both the distances in the real cube example and the Rademacher type \(p\) constant. We define its norm by taking the \(p\)-mean of the real norms of \(\operatorname{Re}(e^{i\theta}z)\) over all rotations. The normalization makes the real subspace isometric, and the \(p\)-mean lets us integrate the real Rademacher inequality with its original constant.

\begin{proof}[Proof of Theorem~\ref{thm:sharp-comparison}]
Areshidze's estimate \cite[p.~5, Corollary~3.1]{aresh} gives the upper bound \(p\). By Lemma~\ref{lem:localized-witness}, finite-dimensional real spaces give ratios \(T_p^E(X_n)/T_p^R(X_n)>\kappa\) for every \(1\le\kappa<p\). For the complex lower bound, equip the complexification of each such real finite-dimensional space with the norm
\[
 \|z\|_{\mathbb C}=
 \left(\frac{\int_0^{2\pi}\|\operatorname{Re}(e^{i\theta}z)\|_X^p\,d\theta}
 {\int_0^{2\pi}|\cos\theta|^p\,d\theta}\right)^{1/p}.
\]
Definiteness follows from continuity in \(\theta\) and evaluation at \(0\) and \(\pi/2\). Complex homogeneity follows by rotating \(\theta\), and the triangle inequality follows from Minkowski's inequality. This norm restricts isometrically to \(X\). Applying the real type inequality to the vectors \(\operatorname{Re}(e^{i\theta}z_j)\) and integrating gives an upper bound of \(T_p^R(X)\) for the complex Rademacher type constant. The real subspace gives the reverse inequality. The real cube witness retains all its distances, so the complex space also violates the comparison with factor \(\kappa\). Letting \(\kappa\uparrow p\) gives the matching lower bound over finite-dimensional spaces in both scalar fields.
\end{proof}

\begingroup
\begin{example}\label{ex:equal-edge-sharpness}
At \(p=2\) the universal factor \(2\) can be approached by cube maps for which every edge has length one and every antipodal pair has the same distance. More precisely, fix \(0<\epsilon<1\). For all sufficiently large \(n\) there is a finite-dimensional real space \(Y\) and a map \(F:\Om_n\to Y\) with
\[
 \|F(x)-F(x^j)\|=1,\qquad
 \|F(x)-F(-x)\|>(2-\epsilon)\,T_2^R(Y)\sqrt n
 \quad (x\in\Om_n,\ j\in[n]).
\]
To see this, take \(X_n,f\) from Lemma~\ref{lem:localized-witness} with \(\kappa=2-\epsilon\). Replace \(f\) by its odd part, which keeps the numerator and can only decrease \(H_2(f)\), and put \(r=H_2(f)/\sqrt n>0\). Let \(\sigma\) cyclically permute the \(n\) coordinates and give \(G=\Om_n\times\{0,\ldots,n-1\}\) uniform probability. In \(Y=L_2(G;X_n)\) set
\[
 F(x)(z,t)=\frac{f(z\,\sigma^t x)}{2r}.
\]
Averaging the squared edge difference over \(z\) averages over all vertices, and averaging over \(t\) gives each coordinate once. Thus its mean is \(4H_2(f)^2/(4nr^2)=1\). Oddness gives \(\|F(x)-F(-x)\|=\|f\|_{L_2(X_n)}/r\) for every \(x\). Applying the Rademacher inequality pointwise and integrating gives \(T_2^R(Y)\le T_2^R(X_n)\), with equality by constant functions. The claimed distance now follows from \eqref{eq:localized-witness}. For instance, \(\epsilon=1/100\) gives \(1.99<T_2^E(Y)/T_2^R(Y)\le2\).
\end{example}
\endgroup

\section{Type inequalities for orthogonal families}\label{sec:orthogonal-tests}
Throughout this section \(X\) is a nonzero real Banach space. The restricted type constants \(B_p\) and \(I_p\) do not increase on passing to a subspace, since each relation is determined by the norm on the span of the two vectors. Write \(\gamma_s=(\E|g|^s)^{1/s}\) for a standard real Gaussian \(g\).

\begin{lemma}\label{lem:unit-gaussian}
Fix \(1<s<2\), and let \(s'=s/(s-1)\). Suppose there is \(a>0\) such that for every \(m\ge1\), \(\delta>0\), and \(z_1,\ldots,z_m\) in the unit ball of \(X\),
\begin{equation}\label{eq:l1-coordinate-lower}
 \left\|\sum_i b_i z_i\right\|\ge\delta\|b\|_1
 \quad\text{for all }b\in\mathbb R^m
 \quad\Longrightarrow\quad
 \delta m^{1/s'}\le a .
\end{equation}
Then for every such family and independent standard Gaussians \(g_i\),
\[
 \E\left\|\sum_{i=1}^m g_i z_i\right\|\le C_s a m^{1/s},
\]
where \(C_s\) depends only on \(s\).
\end{lemma}

\begin{proof}
Put \(E=\operatorname{span}\{z_i\}\) and \(K=\{(\ell(z_i))_{i=1}^m:\ell\in B_{E^*}\}\). For \(t>0\), let \(v_K(t)\) be the maximum cardinality of a coordinate set \(J\) for which the projection of \(K\) onto \(\mathbb R^J\) contains \([-t,t]^J\). Such an inclusion implies \(\|\sum_{i\in J} b_i z_i\|\ge t\|b\|_1\) by duality. By \eqref{eq:l1-coordinate-lower}, \(v_K(t)\le\min\{m,(a/t)^{s'}\}\).

Rudelson and Vershynin's Gaussian entropy bound \cite[Theorem~6.5]{rv} gives
\[
 \E\sup_{w\in K}\sum_{i=1}^m g_i w_i
 \le C\sqrt m\int_0^\infty\sqrt{v_K(t)}\,dt .
\]
The shattering condition in this bound asks for thresholds on a chosen set of coordinates such that every sign pattern is realized with a gap of \(2t\) between its lower and upper levels. For a symmetric convex set, this is equivalent to containing the projected cube above. In fact, the shattering thresholds give a translate of \([-t,t]^J\) in the projection by comparing support functions in every direction. Symmetry and convexity then give the centred cube. Splitting the integral at \(a m^{-1/s'}\) bounds the right-hand side by \(Ca m^{1/s}(1+2/(s'-2))\). Here \(s'>2\), and the supremum on the left equals \(\|\sum_i g_i z_i\|\).
\end{proof}

\begin{lemma}\label{lem:gaussian-blocks}
Let \(E\) be a finite-dimensional real normed space, let \(u_k:\ell_2^{m_k}\to E\), and write \(\sigma_k=\|u_k\|\), \(M_k=\E\|u_k g_k\|\), where \(g_k\) are independent standard Gaussian vectors. Let \(G_k\) be independent \(m_k\times d_k\) standard Gaussian matrices, with \(d_k\ge1\). For \(1\le s\le2\),
\begin{align}
 \E\inf_{t_k\in S^{d_k-1}}\left\|\sum_k u_kG_kt_k\right\|
 &\ge \E\left\|\sum_k u_kg_k\right\|
       -\sum_k \sigma_k(\sqrt{d_k}+\gamma_1),\label{eq:gauss-minmax}\\
 \big(\E\|u_kG_k\|^s\big)^{1/s}
 &\le M_k+\sigma_k(\sqrt{d_k}+1).\label{eq:gauss-operator}
\end{align}
\end{lemma}

\begin{proof}
We use two consequences of Gaussian interpolation. For finite centred Gaussian arrays indexed by \((i,j)\), the expected minimum over \(i\) of the maximum over \(j\) increases if the covariances between distinct rows increase and covariances within each row are fixed. To see this, approximate the maximum in row \(i\) by \(m_i(z)=\alpha^{-1}\log\sum_j e^{\alpha z_{ij}}\) and approximate the minimum by \(H(z)=-\beta^{-1}\log\sum_i e^{-\beta m_i(z)}\). If \(P_i\) and \(Q_{ij}\) are the weights in the two logarithmic averages, the mixed derivative of \(H\) across rows \(i\ne i'\) is \(\beta P_iP_{i'}Q_{ij}Q_{i'j'}\ge0\). Gaussian integration by parts along the interpolation \(\sqrt{s}Z+\sqrt{1-s}W\) gives the derivative of its expectation as one half the sum of covariance differences times the corresponding expected mixed derivatives. The within-row differences vanish. Taking \(\alpha,\beta\to\infty\) proves the comparison. The same interpolation for the maximum alone shows that, at equal variances, the expected maximum decreases when off-diagonal covariances increase, since the mixed derivatives of the logarithmic maximum are nonpositive.

For \eqref{eq:gauss-minmax}, let \(a_k(\ell)=u_k^*\ell\) and \(b_k(\ell)=\|a_k(\ell)\|_2\) for \(\ell\in B_{E^*}\). For \(t=(t_k)\) in the product of the spheres, compare the arrays
\begin{align*}
 Z_{t,\ell}&=\sum_k\langle G_kt_k,a_k(\ell)\rangle
              +\sum_k\eta_k b_k(\ell),\\
 W_{t,\ell}&=\sum_k\langle g_k,a_k(\ell)\rangle
              +\sum_k b_k(\ell)\langle h_k,t_k\rangle ,
\end{align*}
where \(\eta_k\) are standard Gaussians and \(h_k\) are standard Gaussian vectors in \(\mathbb R^{d_k}\). The variables \(\eta_k,h_k,G_k,g_k\) are mutually independent over all blocks. The covariance of \(Z\) minus that of \(W\) at \((t,\ell),(s,\ell')\) equals
\[
 \sum_k (1-\langle t_k,s_k\rangle)
 \big(b_k(\ell)b_k(\ell')-\langle a_k(\ell),a_k(\ell')\rangle\big).
\]
It is nonnegative and vanishes for \(t=s\). The comparison applies first on finite nets and then on the compact index sets by continuity and dominated convergence, using integrable bounds from the finitely many Gaussian norms. The added term in \(Z\) has absolute value at most \(\sum_k\sigma_k|\eta_k|\), and the added term in \(W\) is bounded below by \(-\sum_k \sigma_k\|h_k\|_2\). Since \(\E\|h_k\|_2\le\sqrt{d_k}\), the comparison gives \eqref{eq:gauss-minmax}.

For \eqref{eq:gauss-operator}, omit the index \(k\) and compare the maxima of
\[
 \langle Gt,u^*\ell\rangle+\sigma\eta,
 \qquad
 \langle g,u^*\ell\rangle+\sigma\langle h,t\rangle
\]
over \(t\in S^{d-1}\) and \(\ell\in B_{E^*}\). The variances agree, and the covariance of the first process minus that of the second at \((t,\ell),(v,\ell')\) is \((1-\langle t,v\rangle)(\sigma^2-\langle u^*\ell,u^*\ell'\rangle)\ge0\). The maximum comparison gives \(\E\|uG\|\le M+\sigma\sqrt d\). The function \(G\mapsto\|uG\|\) is \(\sigma\)-Lipschitz in the Hilbert--Schmidt norm, so its variance is at most \(\sigma^2\) by the Gaussian Poincar\'e inequality \cite[p.~180, Proposition~4.1.1]{bakryGentilLedoux} (see also \cite{ledoux}). Indeed, for the Ornstein--Uhlenbeck semigroup \(P_t F(x)=\E F(e^{-t}x+\sqrt{1-e^{-2t}}\,g)\), Gaussian integration by parts and the gradient identity give \(\operatorname{Var}F=2\int_0^\infty\E\|\nabla P_t F\|_2^2\,dt\le\E\|\nabla F\|_2^2\). For Lipschitz \(F\) the bound follows by smooth approximation with the same Lipschitz constant and convergence in Gaussian \(L_2\). Using \(s\le2\) gives \eqref{eq:gauss-operator}.
\end{proof}

\begin{lemma}\label{lem:orthogonal-selection}
Let \(T_j:\mathbb R^{d_j}\to E\), \(1\le j\le k\), be injective linear maps into a finite-dimensional real normed space. If \(d_j\ge2j-1\) for each \(j\), one can select \(t_j\in S^{d_j-1}\) such that \((T_jt_j)_{j=1}^k\) is pairwise mutually Birkhoff--James orthogonal. For pairwise isosceles orthogonality the dimension condition \(d_j\ge j\) suffices.
\end{lemma}

\begin{proof}
For a smooth norm, let \(J(y)\) be the unique norming functional of \(y\ne0\), namely the linear functional satisfying \(\|J(y)\|=1\) and \(J(y)(y)=\|y\|\). This map is continuous and odd. After choosing \(y_i=T_it_i\) for \(i<j\), apply Borsuk--Ulam \cite{borsuk} to the continuous odd map
\[
 t\longmapsto \big(J(y_i)(T_jt),\,J(T_jt)(y_i)\big)_{i<j}
 \quad\text{from }S^{d_j-1}\text{ to }\mathbb R^{2(j-1)}.
\]
The dimension hypothesis ensures a zero, by restriction to an equator if needed. At the resulting \(y_j=T_jt_j\), each earlier norming functional vanishes on \(y_j\), and \(J(y_j)\) vanishes on all earlier vectors. This implies the desired mutual orthogonality.

For an arbitrary norm \(N\), choose smooth norms converging uniformly to it on bounded sets by approximating \(N\) by maxima of finitely many absolute values of norm-one functionals and replacing each maximum by a sufficiently large even-power norm. More explicitly, for any \(\eta>0\) there are \(\ell_1,\ldots,\ell_m\) in the dual unit ball with \((1-\eta)N(y)\le\max_i|\ell_i(y)|\le N(y)\). The functionals span the dual, so \((\sum_i|\ell_i(y)|^{2r})^{1/(2r)}\) is smooth away from zero and converges uniformly to this maximum on bounded sets as \(r\to\infty\). Apply the smooth case to a sequence of such approximations with \(\eta\downarrow0\). Compactness gives a common subsequence along which the selected \(t_j\) converge. Passing to the limit in the defining orthogonality inequalities yields the claim.

For isosceles orthogonality, after choosing \(y_i=T_it_i\) for \(i<j\), take a zero of the continuous odd map
\[
 t\longmapsto\big(\|T_jt+y_i\|-\|T_jt-y_i\|\big)_{i<j}.
\]
This map takes values in \(\mathbb R^{j-1}\), so \(d_j\ge j\) suffices. The same argument works on the unit sphere of any norm on the domain, using its odd radial homeomorphism with the Euclidean sphere.
\end{proof}

\begin{lemma}\label{lem:gaussian-transfer}
Let \(1<s<2\). Suppose that \(\E\|\sum_{i=1}^m g_i z_i\|\le A m^{1/s}\) for every \(z_1,\ldots,z_m\) in the unit ball of \(X\). Then
\[
 T_s^R(X)\le C_s(1+A)B
\]
for either \(B=B_s(X)\) or \(B=I_s(X)\), with \(C_s\) depending only on \(s\).
\end{lemma}

\begin{proof}
Assume \(B<\infty\). Take a finite family \((x_i)\) in \(X\), discard zero vectors, and scale so that \(\sum_i\|x_i\|^s=1\). Let \(E\) be its span. For \(k\ge0\) put
\[
 J_k=\{i:2^{-k-1}<\|x_i\|\le2^{-k}\},\qquad
 d_k=2k+3,\qquad u_k e_i=x_i\quad(i\in J_k),
\]
where \(u_k:\ell_2^{J_k}\to E\) and empty blocks are omitted. As above, let \(\sigma_k=\|u_k\|\) and \(M_k=\E\|u_k g_k\|\). We have
\[
 \sigma_k\le \left(\sum_{i\in J_k}\|x_i\|^2\right)^{1/2}
 \le 2^{-k(1-s/2)}.
\]
In particular, \(D_s=\sum_{k\ge0}2^{-k(1-s/2)}(\sqrt{d_k}+1)<\infty\). Omit each block of rank less than \(d_k\). For such a block \(M_k\le\sigma_k\sqrt{d_k}\) by restricting the Gaussian vector to the orthogonal complement of \(\ker u_k\). The sum of their means is at most \(D_s\).

For retained blocks, \(T_k=u_kG_k\) is almost surely injective. The \(j\)-th retained block in increasing order of \(k\) has \(d_k\ge2j-1\). Put \(L=\inf_{t_k\in S^{d_k-1}}\|\sum_k T_k t_k\|\) over the retained blocks. This is a continuous function of the matrices. By Lemma~\ref{lem:orthogonal-selection}, for each injective sample there are \(y_k=T_k t_k\) forming an admissible family for the chosen constant \(B\). Since each \((\varepsilon_k t_k)\) lies in the same product of spheres,
\[
 L\le \big(\E_\varepsilon\|\textstyle\sum_k\varepsilon_k y_k\|^s\big)^{1/s}
 \le B\left(\sum_k\|T_k\|^s\right)^{1/s}.
\]
The resulting inequality is between measurable functions and uses only the existence of the selected vectors.

By hypothesis \(M_k\le A\,2^{-k}|J_k|^{1/s}\). Since \(\sum_k 2^{-ks}|J_k|\le2^s\), we have \((\sum_k M_k^s)^{1/s}\le2A\). It follows from \eqref{eq:gauss-operator} and Minkowski's inequality on \(\ell_s\) that \(\E L\le B(2A+D_s)\). By \eqref{eq:gauss-minmax}, \(\E\|\sum_{\text{retained }k}u_kg_k\|\le\E L+D_s\). Including the omitted blocks and using \(B\ge1\) gives
\[
 \E\left\|\sum_i g_i x_i\right\|
 \le 2AB+BD_s+2D_s\le(2A+3D_s)B .
\]
If there are no retained blocks, the omitted-block estimate alone gives the same bound.

Finally, for \(F(\varepsilon)=\|\sum_i\varepsilon_i x_i\|\), conditioning on the Gaussian signs gives \(\E F\le\gamma_1^{-1}\E\|\sum_i g_i x_i\|\). The cube Poincar\'e inequality gives \(\operatorname{Var}F\le\sum_i\E|D_iF|^2\le\sum_i\|x_i\|^2\le 1\). The first inequality follows by expanding in Walsh characters, and the second by the reverse triangle inequality. Thus
\[
 (\E F^s)^{1/s}\le(\E F^2)^{1/2}\le\E F+1
 \le \big(1+\gamma_1^{-1}(2A+3D_s)\big)B.
\]
The assertion follows by homogeneity.
\end{proof}

\begin{proof}[Proof of Theorem~\ref{thm:orthogonal-tests}]
Assume \(B_p(X)<\infty\). To check \eqref{eq:l1-coordinate-lower} with \(s=p\), suppose \(z_1,\ldots,z_m\) in the unit ball satisfy \(\|\sum_i b_i z_i\|\ge\delta\|b\|_1\) for all \(b\in\mathbb R^m\), with \(\delta>0\). Put \(N(b)=\|\sum_i b_i z_i\|\). Let \(U\) maximize the absolute determinant among matrices with columns in the unit ball of \(N\). The identity is admissible, so \(|\det U|\ge1\), and each column \(u_j\) has \(N(u_j)=1\). The \(j\)-th coordinate functional of \(U^{-1}\) has norm one for \(N\) by the determinant formula for replacing a column. It annihilates the other columns, so the images of \((u_j)\) in \(X\) are pairwise mutually Birkhoff--James orthogonal. Independent scaling preserves this relation. Conditioning standard Gaussian coordinates on their absolute values gives
\[
 (\E N(Ug)^p)^{1/p}\le B_p(X)\gamma_p m^{1/p}.
\]
If \(r_i\) are the Euclidean lengths of the rows of \(U\), then Hadamard's inequality gives \(\prod_i r_i\ge1\) and hence \(\sum_i r_i\ge m\). Thus \(\E N(Ug)\ge\delta\gamma_1\sum_i r_i\ge\delta\gamma_1 m\). This gives \eqref{eq:l1-coordinate-lower} with \(a=(\gamma_p/\gamma_1)B_p(X)\).

Lemmas~\ref{lem:unit-gaussian} and \ref{lem:gaussian-transfer} now give \(T_p^R(X)\le K_p B_p(X)^2\), using \(B_p(X)\ge1\). The definition gives \(B_p(X)\le T_p^R(X)\) and Section~\ref{sec:introduction} gives \(T_p^R(X)\le T_p^E(X)\le pT_p^R(X)\).
\end{proof}

\begin{proof}[Proof of Theorem~\ref{thm:isosceles-tests}]
Assume \(I_p(X)<\infty\). Take \(z_1,\ldots,z_m\) in the unit ball satisfying the lower estimate in \eqref{eq:l1-coordinate-lower}. Put \(k=\lfloor\sqrt m\rfloor\). Choose disjoint sets of coordinates \(J_j\) with \(|J_j|=j\) for \(1\le j\le k\), which is possible since \(k(k+1)/2\le m\). The lower bound with \(\delta>0\) makes each coordinate map \(b\mapsto\sum_{i\in J_j} b_i z_i\) injective. By the isosceles part of Lemma~\ref{lem:orthogonal-selection} on the \(\ell_1\) unit spheres, there are coefficients \(b^{(j)}\) supported on \(J_j\), each of \(\ell_1\) norm one, such that the vectors \(y_j=\sum_{i\in J_j} b_i^{(j)}z_i\) are pairwise isosceles orthogonal. Each \(\|y_j\|\le1\). Disjoint support gives \(\|\sum_j\varepsilon_j b^{(j)}\|_1=k\) for every sign pattern, so
\[
 \delta k\le
 \big(\E_\varepsilon\|\textstyle\sum_j\varepsilon_j y_j\|^p\big)^{1/p}
 \le I_p(X) k^{1/p}.
\]
Since \(k\ge\sqrt m/2\) and \(1-1/r=(1-1/p)/2\), \eqref{eq:l1-coordinate-lower} follows at \(s=r\) with \(a=2^{1-1/p}I_p(X)\). Lemma~\ref{lem:unit-gaussian} supplies the required Gaussian mean bound. Also \(I_r(X)\le I_p(X)\), by monotonicity of the probability moments and of the sequence norms, applied to each admissible family at its original lengths. As \(1<r<2\), Lemma~\ref{lem:gaussian-transfer} at \(s=r\) and \(B=I_r(X)\) yields the first inequality. The ordinary Enflo upper bound at exponent \(r\) yields the second.
\end{proof}

\begin{remark}\label{rem:orthogonal-edges}
Let \(U_p(X)\) be the least constant in \eqref{eq:rademacher-type} over unit vector families that are pairwise mutually Birkhoff--James and pairwise isosceles orthogonal. These families determine the supremal Rademacher type \(q\) of \(X\). Indeed, the Maurey--Pisier theorem \cite{maurey-pisier} gives copies of finite-dimensional \(\ell_q\) with distortion at most \(2\) in infinite-dimensional \(X\). Select unit vectors in disjoint coordinate blocks of sizes \(3j-2\). At step \(j\), the Borsuk--Ulam theorem applied to the \(3(j-1)\) odd equations \(J(x_i)(x)=J(x)(x_i)=0\) and \(\|x+x_i\|-\|x-x_i\|=0\), \(i<j\), gives the required orthogonality, using smooth approximation as in Lemma~\ref{lem:orthogonal-selection}. Disjoint supports give \(\frac12 n^{1/q}\le U_p(X)n^{1/p}\), so \(U_p(X)<\infty\) implies \(p\le q\). The reverse inequality between the suprema follows from \(U_p(X)\le T_p^R(X)\), and finite-dimensional spaces have type \(2\). In particular, \(I_p(X)<\infty\) implies Enflo type \(s\) for every \(1\le s<p\) \cite{ivv}. Theorem~\ref{thm:isosceles-tests} gives a quantitative bound at \(s=2p/(p+1)\). The endpoint question is whether \(I_p(X)<\infty\) implies Enflo type \(p\). The same implication holds for Pythagorean orthogonality \(\|u-v\|^2=\|u\|^2+\|v\|^2\) \cite{james1945}, even when the test is limited to unit families that also satisfy mutual Birkhoff--James orthogonality. If the supremal type is below \(p\), finite representability gives an \(\ell_q\) model with \(1<q<p\), using the stable representation in \(L_1\) at type \(1\) \cite{maurey-pisier,bdk-stable}. The Poincar\'e--Miranda theorem \cite{kulpa} allows perturbations of \((e_i,c\omega_i)\), with \(c^q=2^{1-q/2}\), satisfying both relations under small norm perturbations. The \(\ell_q\) coordinate block gives growth of order \(n^{1/q}\), contradicting a uniform \(n^{1/p}\) bound.

\textit{(Further observations)}\enspace For \(1<p\le2\), let \(E_p^B(X),E_p^I(X),E_p^{\rm Rob}(X)\) be the constants in \eqref{eq:enflo-type} when the edges at each vertex are pairwise mutually Birkhoff--James, pairwise isosceles, or pairwise Roberts orthogonal, respectively. Roberts orthogonality is the condition \(\|u+tv\|=\|u-tv\|\) for all \(t\in\mathbb R\) \cite{roberts1934}. Let \(A_p^{\rm Rob}(X)\) be the corresponding restricted Rademacher constant. For each pair \((A,E)=(B_p,E_p^B),(I_p,E_p^I),(A_p^{\rm Rob},E_p^{\rm Rob})\), one has \(A(X)\le E(X)\le pA(X)\). The lower bound follows from linear cube maps and the upper bound from \eqref{eq:random-sign-representation}, whose selectors preserve the three relations on edge subfamilies. To see sharpness, use \(f,X_n,c_n\) from Lemma~\ref{lem:localized-witness} with \(1\le\kappa<p\). Put \(r_j=\|D_j f\|_{L_p(X_n)}\) and \(d=1+\lceil\log_2 n\rceil\). Choose distinct \(h_j\in\{-1,1\}^d\) with first coordinate \(1\) and take \(F(x)=((f(xz))_{z\in\Om_n},\sum_j x_jr_jh_j)\) in \(Y=L_p(\Om_n;X_n)\oplus_\infty\ell_\infty^d\). Two distinct \(h_i,h_j\) agree in one coordinate and disagree in another, so every combination of two edges satisfies \(\|aD_iF+bD_jF\|=|a|r_i+|b|r_j\). This gives all three relations and \(H_p(F)=H_p(f)\). The first component preserves the antipodal lower bound. Since \(T_p^R(Y)\le(c_n^p+d)^{1/p}\), each ratio \(E(Y)/A(Y)\) exceeds \(\kappa c_n/(c_n^p+d)^{1/p}\). Letting \(n\to\infty\) and then \(\kappa\uparrow p\) gives sharpness in finite dimensions.
\end{remark}

\section{A limit theorem as the modulus of roundness tends to two}\label{sec:nj-stability}
\begingroup

The classical von Neumann--Jordan constant measures the deviation of a norm from the parallelogram law. Its minimum value is \(1\), attained precisely in inner-product spaces \cite{jordan-von-neumann}. The modified von Neumann--Jordan constant restricts the same quotient to pairs of unit vectors \cite{liu-wang-bao-geometric}. Kato et al.\ extended the classical constant to \(n\)-tuples and related the resulting constants to Rademacher type \cite{kato-takahashi-hashimoto}. To compare the type \(p\) constants, we use the modulus of roundness of Amini-Harandi et al.\ \cite[p.~6, Definition~2.1]{amini-harandi-doust-robertson}, defined by
\begin{equation}\label{eq:nj-p}
 \nu_X(p)=
 \sup_{(x,y)\ne(0,0)}
 \frac{\norm{x+y}^p+\norm{x-y}^p}{\norm{x}^p+\norm{y}^p}.
\end{equation}
The modulus is twice the classical constant at \(p=2\).\footnote{Cui et al.~\cite[p.~2, Definition~1]{cui-huang-hudzik-kaczmarek} divide the quotient in \eqref{eq:nj-p} by \(2^{p-1}\), so their \(C_{\mathrm{NJ}}^{(p)}(X)=2^{1-p}\nu_X(p)\). We use \(\nu_X(p)\) because the sharp four-point factor is \(\nu_X(p)/2\), with a roundness threshold of \(2\) independent of \(p\). In Cui et al.'s convention the threshold is \(2^{2-p}\).} The condition \(\nu_X(p)=2\) is equivalent to \(X\) having roundness \(p\) \cite{amini-harandi-doust-robertson}. For \(1<p\le2\), taking \(y=0\) gives the lower bound \(2\), and the triangle inequality gives the upper bound \(2^p\).

Theorem~\ref{thm:sharp-comparison} gives the sharp universal comparison factor \(p\). Half the modulus in \eqref{eq:nj-p} is the least factor in the four-point roundness inequality, which compares the diagonal and edge sums of the \(p\)-th powers of distances for maps of a two-dimensional cube into \(X\). In the proof below, this factor controls the cost of combining two disjoint sets of cube coordinates. At \(\nu_X(p)=2\), repeated use of the four-point inequality gives \(\Etype(X)=\Rtype(X)=1\) by \eqref{eq:nj-cube}. This makes \(\nu_X(p)-2\) a natural parameter for stability. The next theorem shows that the ratio of the type constants approaches \(1\) uniformly as this parameter tends to zero, independently of the dimension and of the size of \(\Rtype(X)\).

\begin{samepage}
The proof of Theorem~\ref{thm:nj-stability} uses the James constant \(J(X)=\sup_{\|x\|=\|y\|=1}\min\{\|x+y\|,\|x-y\|\}\) \cite{james-uniformly-nonsquare,gao-lau-spheres} and the modulus of convexity of \(W\) on its closed unit ball (for \(0<t\le2\), \(\delta_W(t)=\inf\{1-\|(u+v)/2\|_W:\ u,v\in W,\ \|u\|_W,\|v\|_W\le1,\ \|u-v\|_W\ge t\}\)) \cite[Definition~10.2 and formula~(10.3)]{pisier-martingales}. We use \(J(X)\le\nu_X(p)^{1/p}\) to deduce \(J(X)<2\) from \(\nu_X(p)<2^p\), and \(\delta_W(t)\) to bound the norm decrease under successive coordinate averages. For further applications of these two quantities, see \cite{galego-porto-da-silva,oppenheim-banach-t}.\par
\end{samepage}

\begin{theorem}\label{thm:nj-stability}
Fix \(1<p<2\). Then
\[
 \lim_{c\downarrow2}
 \sup_{\substack{X\ne\{0\}\ \text{real Banach}\\
                 \nu_X(p)\le c,\ \Rtype(X)<\infty}}
 \frac{\Etype(X)}{\Rtype(X)}=1.
\]
\end{theorem}

For a coordinate set \(S\subseteq[n]\), let \(x^S\) be obtained from \(x\) by flipping all coordinates in \(S\), and put \(O_S f(x)=(f(x)-f(x^S))/2\). Thus \(O_{\{i\}}=D_i\) and \(O_{[n]}=P_{\rm odd}\). Write \(A_S f\) for the average of \(f\) over the coordinates in \(S\). For \(0\le u\le1\), let \(P_u\) be the noise operator with Walsh multiplier \(u^{|A|}\) on \(\chi_A\). For the noise operator and its action on Walsh characters, see O'Donnell \cite[p.~55, Definition~2.46 and Proposition~2.47, 2021 edition]{odonnell-abf}. Write \(P_u^S\) for the noise operator acting only on \(S\). All function norms below are Bochner \(L_p\) norms unless an equivalent norm is specified.

\begin{proof}
\begin{samepage}
\emph{Two cube estimates.}
Write \(R=\Rtype(X)\) and \(C=\nu_X(p)/2\). The integrand of \eqref{eq:random-sign-representation} at \(v\) is \(\partial_v P_v f\). Indeed, it has multiplier \(d v^{d-1}\) on Walsh characters of degree \(d>0\) and vanishes on constants. The estimate following \eqref{eq:random-sign-representation}, integrated from \(u\) to \(1\), gives
\begin{equation}\label{eq:nj-noise}
 \norm{f-P_u f}_{L_p(X)}\le pR(1-u)^{1/p}H_p(f).
\end{equation}
We also use \(\Etype(X)\le pR\) from Section~\ref{sec:introduction}.\par
\end{samepage}

The four-point inequality defining roundness \(p\) is stated in \cite[p.~2, Definition~1.1]{amini-harandi-doust-robertson}. To obtain it with a factor \(C\), take four points \(a,b,d,e\in X\) and apply \eqref{eq:nj-p} to \(u+v,u-v\), where \(u=(a-b)/2\) and \(v=(d-e)/2\). Convexity bounds the resulting powers of \(\norm{u+v}\) and \(\norm{u-v}\) by the averages of the corresponding powers of the two cross distances. This yields
\[
 \norm{a-b}^p+\norm{d-e}^p
 \le C\bigl(\norm{a-d}^p+\norm{a-e}^p+\norm{b-d}^p+\norm{b-e}^p\bigr).
\]
The factor \(C\) is optimal, since taking \((a,b,d,e)=(x+y,-x-y,x-y,-x+y)\) recovers half the quotient in \eqref{eq:nj-p}.
For disjoint coordinate sets \(A,B\), use the points
\(f(x),f(x^{A\cup B}),f(x^A),f(x^B)\) in that order and average over \(x\). The two left terms have equal averages, and the cross terms occur in pairs. Thus
\(\norm{O_{A\cup B}f}_{L_p(X)}^p\le C(\norm{O_A f}_{L_p(X)}^p+\norm{O_B f}_{L_p(X)}^p)\).
Recursively partitioning the coordinates into two sets with sizes differing by at most one gives
\begin{equation}\label{eq:nj-cube}
 \norm{O_{[k]}f}_{L_p(X)}\le C^{\lceil\log_2 k\rceil/p}H_p(f)
 \qquad(f:\Om_k\to X).
\end{equation}

\emph{An estimate for random coordinate averages.}
Suppose \(X\) is a subspace of a fixed Banach space \(Z\) with an equivalent uniformly convex norm \(|\cdot|\) satisfying \(\norm{z}\le |z|\le D\norm{z}\).
The space \(W=L_p([0,1];(Z,|\cdot|))\) is uniformly convex by Day's theorem \cite[Theorem~2]{day-uc}. Every cube space with this Bochner norm embeds isometrically into \(W\) using step functions on equal subintervals. Denote the resulting norm by \(\norm{\cdot}_W\) and its modulus of convexity by \(\delta_W\).

For \(k\ge2\), let \(\chi:[N]\to[k]\) be a random coloring with independent uniform values \(\chi(i)\). For \(j\in[k]\), set \(S_j=\chi^{-1}(\{j\})\) and write \(S=S_1\). We denote expectation over \(\chi\) by \(\E_\chi\), using the convention that the subscript specifies the random variable being averaged over \cite[p.~24, Notation~1.4 (2021 edition)]{odonnell-abf}. Take \(1\le L<k\), and put \(u_l=1-l/(k-1)\) for \(0\le l<L\). For every \(h:\Om_N\to X\) and scalars \(|b_l|\le1\), we claim
\begin{equation}\label{eq:nj-random-averages}
 \left(\E_S\norm{(I-A_S)\frac1L\sum_{l=0}^{L-1}b_l P_{u_l}h}_{L_p(X)}^p\right)^{1/p}
 \le \alpha_L\norm{h}_{L_p(X)},
 \qquad \lim_{L\to\infty}\alpha_L=0,
\end{equation}
where \(\alpha_L\) is independent of \(N,k\) and of the subspace \(X\) of \(Z\).

It suffices by homogeneity to consider \(\norm{h}_W=1\). For a fixed coloring, let \(v_l\) be the result of averaging \(h\) over the first \(l\) colors. The next averaging operator maps both \(v_l\) and \(v_{l+1}\) to \(v_{l+1}\) and is contractive. If \(\norm{v_l-v_{l+1}}_W\ge t\), where \(0<t<1\), then \(r=\norm{v_l}_W\) satisfies \(t/2\le r\le1\). Applying uniform convexity to \(v_l/r,v_{l+1}/r\) and then applying the averaging operator gives \(\norm{v_{l+1}}_W\le r(1-\delta_W(t))\). Each such step decreases the norm by at least \(t\delta_W(t)/2\). All steps have length at most \(2\), and the total decrease of the norm is at most \(1\). Consequently
\[
 \frac1L\sum_{l<L}\norm{v_l-v_{l+1}}_W^p
 \le t^p+\frac{2^{p+1}}{Lt\delta_W(t)}.
\]
Given \(S\), let \(T\) be the union of any specified \(l\) of the other color classes. Each coordinate of \(S^c\) belongs to \(T\) independently with probability \(l/(k-1)\). Hence \(P_{u_l}^{S^c}h=\E[A_T h\mid S]\). Since \(P_{u_l}^S\) is contractive and commutes with \(A_S\), Jensen's inequality gives
\[
 \E_S\norm{(I-A_S)P_{u_l}h}_W^p
 \le \E\norm{A_T h-A_{T\cup S}h}_W^p
 =\E\norm{v_l-v_{l+1}}_W^p.
\]
The equality follows by permuting the colors. Applying Jensen's inequality to the average over \(l\) and returning to the original norm proves \eqref{eq:nj-random-averages} with
\[
 \alpha_L=D\inf_{0<t<1}\left(t^p+\frac{2^{p+1}}{Lt\delta_W(t)}\right)^{1/p}.
\]
For each fixed \(t\), the second term tends to zero as \(L\to\infty\). Letting \(t\downarrow0\) then gives \(\alpha_L\to0\).

\emph{The comparison after smoothing.}
It suffices to test odd functions \(f:\Om_n\to X\), since \(P_{\rm odd}\) preserves the antipodal part and contracts the edge norms. Put \(g=L^{-1}\sum_{l<L}P_{u_l}f\). Using \eqref{eq:nj-noise} and contractivity of the noise operators, we obtain
\[
 H_p(g)\le H_p(f),\qquad
 \norm{f-g}_{L_p(X)}\le pR\left(\frac{L-1}{k-1}\right)^{1/p}H_p(f).
\]
Use the same random coloring \(\chi\) and classes \(S_j\) on \([n]\). At a fixed \(x\), multiply all its coordinates in \(S_j\) by the same sign \(\varepsilon_j\). Evaluating \(g\) at the resulting points defines an odd function on \(\Om_k\). Apply \eqref{eq:nj-cube} there and average over \(x\) and the coloring. Multiplication by any fixed set of signs preserves the distribution of \(x\), so
\begin{equation}\label{eq:nj-grouped}
 \norm{g}_{L_p(X)}
 \le C^{\lceil\log_2 k\rceil/p}
 \left(\E_\chi\sum_j\norm{O_{S_j}g}_{L_p(X)}^p\right)^{1/p}.
\end{equation}

For \(i\in S\), write \(h_i(f)(x)=x_iD_i f(x)\), which is independent of \(x_i\), and use the same notation for \(g\). The part of \(g\) that is linear in the signs indexed by \(S\) is
\[
 B_Sg(x)=\sum_{i\in S}x_i A_{S\setminus\{i\}}h_i(g)(x).
\]
Its coefficients depend only on the coordinates outside \(S\). Applying the Rademacher estimate conditionally on those coordinates and using contractivity gives
\(\sum_j\norm{B_{S_j}g}_{L_p(X)}^p\le R^p H_p(g)^p\).
We have \(O_SB_Sg=B_Sg\), and for \(i\in S\),
\(D_i(g-B_Sg)=x_i(I-A_{S\setminus\{i\}})h_i(g)\).
Apply the bound \(\Etype(X)\le pR\) to \(g-B_Sg\) on the cube indexed by \(S\), and integrate over its complementary coordinates. This gives
\[
 \norm{O_Sg-B_Sg}_{L_p(X)}^p
 \le (pR)^p\sum_{i\in S}\norm{(I-A_{S\setminus\{i\}})h_i(g)}_{L_p(X)}^p.
\]
For an empty color class both sides are zero. The definition of \(P_u\) gives
\(h_i(g)=L^{-1}\sum_{l<L}u_l P_{u_l}^{[n]\setminus\{i\}} h_i(f)\).
Conditioned on the color of \(i\), the remaining coordinates of its color class have the law of \(S\) in \eqref{eq:nj-random-averages} on \(n-1\) coordinates. Applying that estimate with \(b_l=u_l\), and summing over \(i\), bounds
\(\E_\chi\sum_j\norm{O_{S_j}g-B_{S_j}g}_{L_p(X)}^p\) by \((pR\alpha_L)^p H_p(f)^p\).
Minkowski's inequality in \eqref{eq:nj-grouped}, followed by the bound on \(f-g\), now gives
\begin{equation}\label{eq:nj-stability-estimate}
 \frac{\Etype(X)}{\Rtype(X)}
 \le p\left(\frac{L-1}{k-1}\right)^{1/p}
 + C^{\lceil\log_2 k\rceil/p}(1+p\alpha_L).
\end{equation}

\emph{Passage to the limit.}
Suppose there are \(\varepsilon>0\) and nonzero real Banach spaces \(X_m\) with \(\nu_{X_m}(p)\to2\), finite \(\Rtype(X_m)\), and \(\Etype(X_m)/\Rtype(X_m)\ge1+\varepsilon\).
Put \(C_m=\nu_{X_m}(p)/2\). After discarding finitely many terms, \(C_m\le c_0<2^{p-1}\). Let \(Z=(\bigoplus_m X_m)_{\ell_p}\) with the original norms on the summands. For \(z=(z_m)\) and \(w=(w_m)\) in \(Z\), summing the two-vector bound in the summands gives
\[
 \|z+w\|_Z^p+\|z-w\|_Z^p
 =\sum_m\bigl(\|z_m+w_m\|_{X_m}^p+\|z_m-w_m\|_{X_m}^p\bigr)
 \le 2c_0\bigl(\|z\|_Z^p+\|w\|_Z^p\bigr).
\]
In particular, for unit vectors \(z,w\) in \(Z\) the smaller of \(\|z+w\|_Z^p\) and \(\|z-w\|_Z^p\) is at most half their sum, hence at most \(2c_0\). Taking the supremum over all such pairs gives \(J(Z)\le(2c_0)^{1/p}<2\). Thus every unit pair has at least one of \(\|z+w\|_Z\) and \(\|z-w\|_Z\) bounded by \(2-\eta\), with the same \(\eta=2-(2c_0)^{1/p}>0\). This is uniform nonsquareness, so \(Z\) admits an equivalent uniformly convex norm \cite[Theorem~11.37]{pisier-martingales}. The same sequence \((\alpha_L)\) in \eqref{eq:nj-stability-estimate} applies to every summand. All type constants in that estimate use the original norms. For fixed \(L<k\), first take the upper limit as \(m\to\infty\) in \eqref{eq:nj-stability-estimate}. We obtain
\[
 \limsup_{m\to\infty}\frac{\Etype(X_m)}{\Rtype(X_m)}
 \le p\left(\frac{L-1}{k-1}\right)^{1/p}+1+p\alpha_L.
\]
Letting \(k\to\infty\) and then \(L\to\infty\) bounds the upper limit by \(1\), a contradiction. Linear cube maps give \(\Rtype(X)\le\Etype(X)\). The real line has \(\nu_{\mathbb R}(p)=2\) by the scalar Clarkson inequality, and \eqref{eq:nj-cube} gives both type constants equal to \(1\) there. The asserted limit follows.
\end{proof}

For the Birkhoff--James restriction, the smoothing argument in Theorem~\ref{thm:nj-stability} gives the following estimate for admissible \(f\), with \(S\) one of \(k\) random color classes and \(0<u<1\) fixed:
\begin{equation}\label{eq:nj-b-smoothed}
 k^{1/p}\|O_SP_u f\|_{L_p(S,x)}\le (1+o_k(1)) B_p(X_m) H_p(f)
\end{equation}
along any sequence \(X_m\) with \(\nu_{X_m}(p)\to2\). Put \(B=B_p(X_m)\) and \(C_m=\nu_{X_m}(p)/2\). For the Walsh projection \(L_S\) onto degree one in \(S\), average \(\sum_{i\in S}x_i y_iD_i f(y_S,x_{S^c})\) over \(y_S\). Jensen's inequality lets us apply the \(B\)-estimate at a single vertex before averaging, giving the leading term \(BH_p(f)\) after scaling by \(k^{1/p}\), also for \(L_S P_u f=P_u L_S f\). For \(Q_S=O_S-L_S\), a fixed degree-\(d\) Walsh projection is \(d^d/d!\) times the average of products of \(d\) first-degree projections on random groups. These products have norm at most \(K^d\) by Pisier's theorem \cite{pisier-holomorphic} on the common superreflexive space above. Each surviving pattern meets at least three groups. In one hit group the slice estimate of Remark~\ref{rem:orthogonal-edges} gives a factor \(pB k^{-1/p}\) before projection. In a second group, the common Bochner space has cotype \(q<\infty\) \cite{maurey-pisier}. This means \((\sum_b\|F_b\|^q)^{1/q}\) is bounded by a constant times \((\E_\varepsilon\|\sum_b\varepsilon_b F_b\|^2)^{1/2}\) for all finite families \(F_b\) in that space. Apply this to the sums indexed by color. Sign invariance and exchangeability then give \(O(k^{-1/q})\) in the \(L_p\) average upon retaining one color. The other block projections are contractions. This proves the scaled decay on each fixed degree uniformly in \(m,n\). It extends to \(P_u\) by dimension-uniform approximation in the span of the degree projections, from the convergent Walsh expansion near \(u=0\). Indeed, the semigroup \(P_{e^{-t}}\) extends holomorphically in \(t\) to a common complex sector, with operator bounds uniform in the cube dimension \cite{pisier-holomorphic}. The identity theorem in the quotient by the closed span of the degree projections then extends the approximation to all \(0<u<1\). The approximants commute with \(Q_S\), and the slice bounds give \(k^{1/p}\|Q_Sf\|_{L_p(S,x)}\le(p+1)BH_p(f)\).

\begin{corollary}\label{thm:nj-b-stability}
Fix \(1<p<2\). With \(X\) ranging over nonzero real Banach spaces,
\[
 \lim_{c\downarrow2}
 \sup_{\substack{\nu_X(p)\le c\\B_p(X)<\infty}}
 \frac{E_p^B(X)}{B_p(X)}=1.
\]
\end{corollary}

\begin{proof}
Use \eqref{eq:nj-b-smoothed} in the grouped cube estimate \eqref{eq:nj-grouped} for \(P_uf\). For \(f-P_uf\), integrate \eqref{eq:random-sign-representation} from \(u\) to \(1\) with the restricted type estimate. Together they give
\[
 \frac{E_p^B(X_m)}{B_p(X_m)}
 \le p(1-u)^{1/p}+C_m^{\lceil\log_2 k\rceil/p}(1+o_k(1)).
\]
Take the upper limit as \(m\to\infty\), then let \(k\to\infty\) and \(u\uparrow1\). Linear cube maps give the reverse bound.
\end{proof}

\begin{remark}\label{cor:nj-b-edges}
Fix \(1<p<2\). Over nonzero real Banach spaces,
\[
 \lim_{c\downarrow2}
 \sup_{\substack{\nu_X(p)\le c\\T_p^R(X)<\infty}}
 \frac{E_p^B(X)}{T_p^R(X)}=1.
\]
The upper bound follows from \(E_p^B(X)\le T_p^E(X)\) and Theorem~\ref{thm:nj-stability}, and a one-coordinate cube in \(\mathbb R\) attains \(1\).
\end{remark}

\endgroup

The modulus of roundness gives a geometric answer to the stability question for the type constants. As \(\nu_X(p)\) tends to \(2\), the Enflo to Rademacher type ratio tends uniformly to \(1\), both for the ordinary constants and for their restrictions to mutual Birkhoff--James orthogonality. The proof uses the James constant to obtain uniform nonsquareness on a common direct sum and the modulus of convexity to control the error from successive coordinate averages. A quantitative question is to determine the optimal rate in terms of \(\nu_X(p)-2\). We expect that geometric constants can also help establish stability results for other comparisons between linear and metric inequalities.

\section*{Acknowledgements}
Thanks to all the members of the Functional Analysis Research team of the College of Mathematics and Statistics of Anqing Normal University for their discussion and correction of the difficulties and errors encountered in this paper. 

\end{document}